\documentclass[12pt]{article}

\usepackage{amsfonts}
\usepackage{mathtools}
\usepackage{amsmath}
\usepackage{amssymb}
\usepackage{amsthm}
\usepackage{enumerate}
\usepackage{graphicx}

\def \Z {\mathbb Z}

\def\P{\mathbb P}

\def\ep{\varepsilon}

\def\Ph1{P^{(h_1)}}
\def\Ph2{P^{(h_2)}}

\def\b0{{\bf 0}}

\newtheorem{thm}{Theorem}[section]
\newtheorem{prop}[thm]{Proposition}
\newtheorem{lem}[thm]{Lemma}

\theoremstyle{plain}

\newtheorem{defn}[thm]{Definition}

\newtheorem{obs}[thm]{Observation}

\theoremstyle{definition}

\begin{document}

\title{An upper bound on the two-arms exponent for critical percolation on $\Z^d$}

\author{J. van den Berg\footnotemark[2] 
%\footnote{During part of this project JvdB was affiliated
%as visiting professor at NYU Abu Dhabi} 
$\,$ and D.G.P. van Engelenburg\footnotemark[4] \footnote{This paper arose from a part of the master
thesis work at the University of Amsterdam and the VU University Amsterdam by DGPvE
in the period September 2019 - February 2020, supervised by JvdB.}\\
%{\it\bf Draft, Friday December 20, 2019}
\small \\
  \small \footnotemark[2] CWI and VU University Amsterdam; email: J.van.den.Berg@cwi.nl \\
  \small \footnotemark[4] University of Vienna; email: diederik.van.engelenburg@univie.ac.at
  }
\date{}
\maketitle

\begin{abstract}
Consider critical site percolation on $\Z^d$ with $d \geq 2$.
Cerf (2015) pointed out that from classical work by Aizenman, Kesten and Newman (1987) 
and Gandolfi, Grimmett and Russo (1988) one can obtain that the two-arms exponent is at least $1/2$.
The paper by Cerf slightly improves that lower bound.

 Except for $d=2$ and for high $d$, no {\it upper} bound for this
exponent seems to be known in the literature so far (not even implicity).
We show that the distance-$n$ two-arms probability is at least $c \, n^{-(d^2 + 4 d -2)}$ (with $c >0$ a constant
which depends on $d$), thus giving an upper bound $d^2 + 4 d -2$ for the above mentioned exponent.

\end{abstract}
{\it Key words and phrases:} Critical percolation, critical exponent. \\
{\it AMS 2000 subject classifications.} Primary 60K35; secondary 82B43.

\newpage
\begin{section}{Introduction and main results}

Consider nearest-neighbour site percolation with parameter $p$ on $\Z^d$, with $d \geq 2$.
Each site is, independently of the others, open with probability $p$
and closed with probability
$1-p$. The corresponding product measure is denoted by $P_p$. We use the standard notation $\theta(p)$ for the probability that a given site belongs to an infinite
open cluster, and   $p_c = p_c(d):= \sup\{p \in [0,1] \, : \, \theta(p) = 0\}$ for the critical probability.

Let $\Lambda(n)$ be the box $[-n,n]^d$, and let $\partial \Lambda(n)$ be its boundary $\{x \in \Lambda(n) \, : \,
\|x\| = n\}$, where $\|x\| := \max_{1\leq i \leq d} |x_i|$.
Following \cite{C15}, we use the notation ``two-arms$(0,n)$" for the event that in $\Lambda(n)$
there are two distinct open clusters connecting neighbors of $0$ to 
$\partial \Lambda(n)$. More precisely, it is the event that there are open paths $\pi$ and $\pi'$ with the following properties:
$\pi$ connects a neighbor of $0$ with $\partial \Lambda(n)$; $\pi'$ connects a neighbor of $0$ with $\partial \Lambda(n)$;
and there is no open path in $\Lambda(n)$ connecting $\pi$ with $\pi'$.

\smallskip
In two dimensions, two-arms events are, somewhat confusing, in fact  events involving {\it four arms} with alternating `colour', 
and have been widely studied since Kesten's
classical work \cite{K80}.

In general dimensions two-arms events are related to the uniqueness of the infinite cluster. Indeed, it is straightforward that
\begin{eqnarray} \nonumber
 \, & & \lim_{n \rightarrow \infty} P_p(\text{ two-arms }(0,n)) \\ \nonumber
 &=&  P_p(\exists \text{ two neighbours of } 0 \text{ that belong to distinct infinite open clusters}),
\end{eqnarray}
which is known to be $0$ (`uniqueness of the infinite cluster').
The paper \cite{BK89} by Burton and Keane gives the most elegant proof of uniqueness of the infinite cluster, but does not give a bound
for the speed of 
convergence above. Cerf \cite{C15} pointed out that from earlier uniqueness proofs (in \cite{AKN87} and \cite{GGR88})
%yield a `power-law' upper bound on two-arms probabilities,
%in particular 
one can obtain that in any
dimension the two-arms exponent is at least $1/2$. One of Cerf's main results in the above mentioned paper
is a clever computation which slightly improves that bound.
More precisely, he showed that
$$\limsup_{n \rightarrow \infty} \frac{\log P_{p_c}(\text{ two-arms }(0,n))}{\log n} \leq - \,  \frac{2d^2 + 3d -3}{4d^2 + 5d -5},$$
i.e. the two-arms exponent is at least $(2d^2 + 3d -3) / (4d^2 + 5d -5)$.

 To our knowledge, the literature, so far, does not provide {\it upper} bounds on the two-arms exponent,
except for $d=2$ and for `high' $d$.
Our main result, Theorem \ref{Mainthm} below, gives the upper bound $d^2 + 4d -2$.
Somewhat remarkably, one of the main ingredients of the proof of our theorem is an intermediate result by Cerf which he used for
his above mentioned lower bound for the exponent.

Our upper bound on the exponent (as well as Cerf's lower bound) is essentially only relevant for `intermediate' dimensions: 
For site percolation on the triangular lattice the exact value $5/4$ has been proved
for the exponent (see \cite{SW01}), and this is also believed to hold for other nice $2D$ lattices, including the square lattice.
In \cite{KN11} it has been proved that for `high' dimensions 
the exponent is $4$, and this is believed to hold for all $d \geq 7$.

Two-arms events are also (at least potentially) important for the study of one of the
main conjectures in percolation theory, namely that $\theta(p_c) =0$
(this conjecture has been proved for $d=2$ and for high $d$ only). See the last paragraph of Section 1 in \cite{C15}.

Some other examples where results on two-arms probabilities are important, are works about upper bounds for the percolation correlation
length near criticality (\cite{DKT20})
and about the convergence rate in a central limit theorem for minimal spanning trees \cite{CS17}.

\medskip
We now state our main result more precisely, where (here and throughout the rest of the paper) 
{\it we use the shorthand notation $\P$ for $P_{p_c}$}.

\begin{thm}\label{Mainthm}
There is a constant $c>0$, which depends only on the dimension $d$, such that for all $n \geq 1$,

\begin{equation}\label{eq-Mainthm}
\P(\text{\em{two-arms}}(0,n)) \geq c \, n^{-(d^2 + 4d -2)}. 
\end{equation}
\end{thm}

A crucial part in the proof of Theorem \ref{Mainthm} is the following proposition, which may also be interesting in
itself.
We first introduce more notation. 
For $m \leq n$, we write $A_2(m,n)$ for the event that in $\Lambda(n)$ there
are two distinct open clusters connecting $\Lambda(m)$ with $\partial \Lambda(n)$. \\
If $W \subset \Z^d$ and $U , V \subset W$, we use the notation $\{U \leftrightarrow V\}$ for the event that
there is an open path from $U$ to $V$, and the notation $\{U \leftrightarrow V \text{ in } W\}$ for the event that there is an open path from
$U$ to $V$ of which every site is in $W$.

\begin{prop}\label{MainProp}
%Let $d\geq 2$. Let $M \geq 2$ and $c, \tilde c$ be such that, for all $n$,
%
%\begin{equation}\label{eq-prop-cond}
%\P(A_2(n, Mn)) \leq c \, n^{\tilde c} \, \P(\text{\em{two-arms}}(0,n)).
%\end{equation}
%Then there is a $\hat c > 0$ such that for all $n$,
%
%\begin{equation}\label{eq-prop-concl}
%\P(\text{\em{two-arms}}(0,n)) \geq \hat c n^{- \tilde c}.
%\end{equation}
%\end{prop}

For all $M \geq 2$ there is a $\delta = \delta(d, M) > 0$ and a sequence of integers
$1 = n_1 < n_2 < \cdots$ such that $\max_i \frac{n_{i+1}}{n_i} \leq 8 M$, and

\begin{equation}\label{eq-prop-opt1}
\P(A_2(n_i, M n_i)) \geq \delta, \,\, \text{ for all } i.
\end{equation}

%\begin{equation}
%\max\left(\P\left(A_2(4n, 4Mn)\right), \P\left(A_2(8 M n, 8 M^2 n)\right)\right) \geq \delta.
%\end{equation}

\end{prop}

We postpone the proof of Proposition \ref{MainProp} to Section \ref{Proof-Prop}, and now show how Theorem \ref{Mainthm}
follows from that proposition, combined with (a straightforward generalization of) the following result by Cerf and a lower
bound for `restricted' one-arm probabilities.

\begin{lem}\label{l-cerf7.2-7.3} {\em\small\bf(Case $l = 2n$ and $p=p_c$ of Corollary 7.3 in \cite{C15})} \\
Let $d \geq 2$. There is a $C = C(d) > 0$ such that for all $n \geq 1$,
\begin{equation}\label{eq-c-7.3}
\P(A_2(n, 3 n)) \leq
\frac{C n^{4d-2} \P(\text{\em{two-arms}}(0,n))}{\min_{x, y \in \partial\Lambda(n)} \P(x \leftrightarrow y \text{ in } \Lambda(2n))}.
\end{equation}
\end{lem}

\medskip\noindent
{\it Proof of Theorem \ref{Mainthm} from Proposition \ref{MainProp}}. \\
We first prove a slightly weaker version (with a larger exponent) and then point out how to obtain the actual theorem.
Choose $M = 3$. Let $n \geq 1$. Take $i$ such that $n \leq n_i \leq 24 n (= 8 M n)$.
By 
\eqref{eq-prop-opt1} and \eqref{eq-c-7.3} we have, 

$$\delta \leq \P(A_2(n_i, 3 n_i)) \leq 
\frac{C \, (n_i)^{4d-2} \, \P(\text{two-arms}(0,n_i))}{\min_{x, y \in \partial\Lambda( n_i)} \P(x \leftrightarrow y \text{ in }
\Lambda( 2 n_i))}.$$
By Lemma 6.1 in \cite{C15}, the denominator in the expression above is at least some constant (which depends on $d$ only) 
times $(n_i)^{-2 d (d-1)}$.
Hence

\begin{eqnarray}\label{eq-cons-opt1} \nonumber
\P(\text{two-arms}(0,n)) &\geq& \P(\text{two-arms}(0,n_i)) \geq C' \delta (n_i)^{-((4d-2) + 2d(d-1))} \\
&\geq& C'' \, n^{-((4d-2) + 2d(d-1))},
\end{eqnarray}
where $C'$ and $C''$ are positive constants wich depend only on $d$.
%Alternatively, if 
%\eqref{eq-prop-opt2} holds,
%we get, (similarly as above, again by using \eqref{eq-c-7.3}),
%
%$$\delta \leq \P(A_2(24 n, 72 n)) \leq
%\frac{C \, (24 n)^{4d-2} \,  \P(\text{two-arms}(0,24 n))}{\min_{x, y \in \partial\Lambda(24 n)} \P(x \leftrightarrow y \in \Lambda(48 n))}.$$
%
%Using again the above mentioned Lemma 6.1 in \cite{C15}, this gives a lower bound for $\P(\text{two-arms}(0,n))$ with the same exponent as 
%\eqref{eq-cons-opt1} but with a different multiplicative constant.
This proves a slightly weaker form of Theorem \ref{Mainthm}: with exponent $(4d-2) + 2d(d-1) = 2d^2 + 2d -2$ instead of $d^2 + 4d -2$.

The exponent in Theorem 1 can be obtained by using Lemma \ref{l-g-cerf7.2-7.3} below
(a generalization of Lemma \ref{l-cerf7.2-7.3} with practically
the same proof), and Lemma \ref{lem-BD} below (which is an improvement in \cite{BD20} of Cerf's Lemma 6.1 mentioned above).

\begin{lem} \label{l-g-cerf7.2-7.3}
For all $M \geq 2$ there is a $c_1 > 0$ such that for all $n \geq 1$,
\begin{equation} \label{eq-g-cerf7.2-7.3}
\P(A_2(n, M n)) \leq
\frac{c_1 \, n^{4d-2} \, \P(\text{\em{two-arms}}(0,n))}{\min_{x, y \in \partial\Lambda(n)} \P(x \leftrightarrow y \text{\em{ in }} \Lambda((M-1) n))}.
\end{equation}
\end{lem}

\begin{lem} \label{lem-BD} {\em\bf\small(Corollary 1.2 in \cite{BD20})}
There is a constant $c_2 = c_2(d) > 0$ such that
for all $n \geq 1$, 
\begin{equation}\label{eq-BD}
\min_{x,y \in \partial \Lambda(n)} \P(x \leftrightarrow y \text{\em{ in} } \Lambda(9 n)) \geq c_2 n^{-d^2}.
\end{equation}
\end{lem}

\noindent
By now choosing $M=10$, and repeating the same argument which led to the earlier mentioned weaker version of Theorem \ref{Mainthm} (but now using
Lemma's \ref{l-g-cerf7.2-7.3} and \ref{lem-BD} instead of Lemma \ref{l-cerf7.2-7.3} and Cerf's Lemma 6.1, respectively) we get, for some constant $c>0$, 

$$\P(\text{two-arms}(0,n)) \geq c\, n^{-((4d-2) + d^2)} = c \, n^{-(d^2+4d-2)},$$
which completes the proof of Theorem \ref{Mainthm}.
\qed

\medskip\noindent
{\bf Remark:} From the above it is clear that there are potentially at least two strategies to improve the exponent in
Theorem \ref{Mainthm}: By improving the exponent $4d-2$ in \eqref{eq-c-7.3} (and in its generalization \eqref{eq-g-cerf7.2-7.3});
or by further improving Lemma \ref{lem-BD} (which
itself, as mentioned before, was an improvement of Lemma 6.1 in \cite{C15}).

\end{section}

\begin{section}{Proof of Proposition \ref{MainProp}} \label{Proof-Prop}

We will use the following lemma, where we use the notation $V(k,m)$ for the event that there is an open crossing in the $d$-direction of the rectangle
$[0,m]^{d-1} \times [0,k]$. (Typically we will use this with $k < m$, so that the crossing is in the `easy' direction). \\
Subsets of $\Z^d$ will sometimes be denoted in a somewhat informal way; for instance we will simply write $\{z_d = 0\}$ for the
set $\{z \in \Z^d \, : \, z_d = 0\}$.

\begin{lem}\label{l-two-opts}
For all $M \geq 2$ there is an $\ep = \ep(M,d) >0$ such that, for all $n \geq2$ the following holds:

\begin{equation} \label{l-option1}
\P(A_2(4n, 4Mn)) \geq \ep,
\end{equation}
or
\begin{equation}\label{l-option2}
\P(V(16Mn, 16M^2n)) \leq 1 - \ep.
\end{equation}
\end{lem}

\medskip\noindent
{\em Proof of Lemma \ref{l-two-opts}} \\
The proof is by contradiction. Suppose the lemma is not true. Then there is an $M \geq 2$ such that for all $\ep >0$ there is an $n \geq 2$ with

\begin{equation}\label{eq-l-opt1}
\P(A_2(4n, 4Mn)) < \ep,
\end{equation}
and

\begin{equation}\label{eq-l-opt2}
\P(V(16Mn, 16M^2 n)) > 1 - \ep.
\end{equation}

Fix such and $M$ and choose $\ep$ very small (how small will become clear later). We will show that we get a contradiction. Take $n$ such that 
\eqref{eq-l-opt1} and \eqref{eq-l-opt2} hold.
First, observe that the event $V(16Mn, 16M^2 n))$ is contained in the event that the `face' $[0,16M^2 n]^{d-1} \times \{0\}$ is connected
to the hyperplane $\{z_d = 16 M n\}$ by an open path in the half-space $\{z_d \geq 0\}$. 
By partitioning the above mentioned `face' in `subfaces' of length $2n$, and using a version of the standard `square-root trick' (see e.g. \cite{G99}), we get
that the probability
of having an open path in the half-space $\{z_d \geq 0\}$ from $[0,2n]^{d-1} \times \{0\}$ to the hyperplane 
$\{z_d = 16 M n\}$ is `large', and hence that 
$\P(\Lambda(n) \leftrightarrow \partial\Lambda(16 M n))$ is `large.'
More precisely,

\begin{equation}\label{eq-g1}
\P(\Lambda(n) \leftrightarrow \partial\Lambda(16 M n)) \geq f(1 - \ep),
\end{equation}
where $f \, : \, [0,1] \rightarrow [0,1]$ is an increasing function with the property that $f(a)$ converges to $1$ as $a \rightarrow 1$. \\
{\bf Remark}: $f$ depends on $M$ but recall that we have fixed $M$ earlier in this proof.

\smallskip
We proceed the proof with a renormalization argument.
For $x \in \Z^d$, we use the notation $\Lambda(x; n)$ for the box obtained by translating 
$\Lambda(n)$ over the vector $x$.
Note that if $x, y \in \Z^d$ are neighbours, then

\begin{equation}\label{eq-obs-1}
\Lambda(n x; 4 n) \supset \Lambda(n y;n),
\end{equation}
and

\begin{equation}\label{eq-obs-2}
\Lambda(n y; 16 M n) \supset \Lambda(n x;4M n).
\end{equation}

Let, for each $x \in \Z^d$, $W(x)$ be the random set of vertices 
$\{z \in \Lambda(n x; n)\, : \, z \leftrightarrow \partial\Lambda(n x; 16 M n)\}$.

%By \eqref{eq-g1} we have that
%
%\begin{equation}\label{eq-g1-bis}
%\P(W(x) \neq \emptyset) \geq f(1 - \ep).
%\end{equation}
%Recall that this is close to $1$ if $\ep$ is close to $0$.

\begin{defn}\label{def-good}
Call a point $x \in \Z^d$ {\em good} if the following holds:

\begin{equation}\label{eq-defgood1}
W(x) \neq \emptyset
\end{equation}
and any two open paths connecting
$\Lambda(nx; 4n)$ with $\partial\Lambda(nx; 4 M n)$ are connected by an open path in $\Lambda(n x; 4 M n)$.
\end{defn}

So, for each site $x$, 
\begin{eqnarray}\label{eq-good-bnd}
\P(x \text{ is good}) & = & \P(0 \text{ is good}) \\ \nonumber
 &\geq& \P\left(\Lambda(n) \leftrightarrow \partial\Lambda(16 M n)\right) - \P\left(A_2(4n, 4Mn)\right) \\ \nonumber
 &\geq& f(1 - \ep) - \ep,
\end{eqnarray}
where the last inequality comes from \eqref{eq-l-opt1} and \eqref{eq-g1}.
Recall that the last expression converges to $1$ as $\ep \rightarrow 0$.

Further, suppose that two neighbours $x$ and $y$ are both good. Then, by definition, $W(x) \neq \emptyset$ and
$W(y) \neq \emptyset$. Moreover, if $v \in W(x)$ and $z \in W(y)$ we have (again by definition), that
$v \leftrightarrow \partial\Lambda(n x; 16 M n)$ and $z \leftrightarrow \partial\Lambda(n y; 16 M n)$.
The corresponding path for $v$ obviously crosses the annulus between $\Lambda(nx;4n)$ and $\partial\Lambda(nx;4M n)$.
By \eqref{eq-obs-1} and \eqref{eq-obs-2} the corresponding path for $z$ also crosses that annulus.
Hence by the second property in the definition of {\it good} we have that $v$ and $z$ are connected by an open path.
Concluding, if two neighbours $x$ and $y$ are both good, then $W(x)$ and $W(y)$ are non-empty and every point of $W(x)$ has an open path
to every point of $W(y)$. This gives immediately the following:

\begin{obs}\label{obs-inf}
If there is an infinite path of good points, then there is an infinite open cluster.
\end{obs}

The remaining part of the proof of Lemma \ref{l-two-opts} is quite standard and we summarize it as follows:
By taking $\ep$ sufficiently small, we can make $\P(x \text{ is good })$ as large as we want (by \eqref{eq-good-bnd}).
Moreover, the process $(x \text{ is good } \, : \, x \in \Z^d)$ is finite-dependent (and the `dependence range' depends only on $M$, which we have
fixed). Hence, using general domination results in \cite{LSS97}, we can take $\ep$ so small that the above mentioned process
dominates a supercritical Bernoulli site percolation process. Finally, again using that whether or not $x$ is {\it good} depends only
on the states (open/closed) of the vertices in $\Lambda(nx;16 Mn)$, the function $p \rightarrow P_p(x \text{ is good})$ is continuous, 
and hence the above domination
also holds for some $p$ smaller than $p_c$. So, by Observation \ref{obs-inf}, $\theta(p) > 0$ for some $p < p_c$.

This is the desired contradiction and completes the proof of Lemma \ref{l-two-opts}. \qed

\begin{subsection}{Completion of the proof of Proposition \ref{MainProp}}
Fix an $M \geq 2$. Take $\ep = \ep(M)$ as in Lemma \ref{l-two-opts}. Let $n \geq 2$.
Suppose that \eqref{l-option2} holds. Then (as is easy to see by a standard argument) there is an $\ep^* = \ep^*(\ep) > 0$ such that the
event $V(16 M n -2, 16 M^2 n)$ has probability at most $1 - \ep^*$.  (We need this tiny modification of the event $V(16 M n, 16 M^2 n)$
to get independence of the events $E_1$, $F_2$ and $F_3$ defined below).
Hence with probability $\geq \ep^*$ there is a `barrier' of closed sites in the box
$[-8 M^2 n,8 M^2 n]^{d-1} \times [-8 Mn, 8 M n]$ which separates
top from bottom of that box (and which itself does not touch that top and bottom). Call this event $E_1$.

See Figure 1. Let $E_2$ be the event that there is an open crossing in the $d$-direction of the box 
$[-8 M^2 n,8 M^2 n]^{d-1} \times [8 Mn, 8 M^2 n]$, and $E_3$ the analog of this event for the box obtained
by refelection in the hyperplane $\{z_d = 0\}$. 
Note that the annulus between $\Lambda(8 M n)$ and $\partial\Lambda(8 M^2 n)$
is the union of $2 d$ `copies' (obtained by rotations and translations)
of the box in the definition of $E_2$, and that each path from $\Lambda(8 M n)$ to $\partial \Lambda(8 M^2 n)$ crosses at least one of those
boxes in the `easy' direction. Further, by a classical argument going back to the work by Kesten \cite{K81},
the probability that there is such an open path is larger than
some constant $\tilde\ep = \tilde\ep(M,d)>0$. Hence, 
$\P(E_2) = \P(E_3) \geq \tilde\ep / (2 d)$.
%From this, with a similar argument (version of the `square-root trick') as earlier in this paper, we get that,
By partitioning $[-8 M^2 n,8 M^2 n]^{d-1} \times \{8 Mn\}$ in translates of $[-8M n, 8 M n]^{d-1} \times \{8 Mn\}$, we get that
there is a constant $\ep' = \ep'(M,d)  > 0$ such that 
$\P(F_2) = \P(F_3) > \ep'$,
where the event $F_2$ is defined by
$$F_2 = \left\{\begin{gathered}
[-8M n, 8 M n]^{d-1} \times \{8 M n\} \leftrightarrow \Z^{d-1} \times \{8 M^2 n\} \\
\text{ in the half-space } \{z_d \geq 8 M n\}
\end{gathered}\right\}$$
%$$F_2 = \left\{[-8M n, 8 M n]^{d-1} \times \{8 M n\} \leftrightarrow \Z^{d-1} \times \{8 M^2 n\} \text{ in the half-space } \{z_d \geq 8 M n\} \right\},$$
and $F_3$ is the event obtained from $F_2$ by reflection in the  hyperplane $\{z_d = 0\}$.

%$$\P([-8M n, 8 M n]^{d-1} \times \{8 M n\} \leftrightarrow \Z^{d-1} \times \{8 M^2 n\} \text{ in the `upper' half-space }) > \ep',$$
%and hence that 
%where $F_2$ is the event that there is an open path in the box $[-8 M^2 n,8 M^2 n]^{d-1} \times [8 Mn, 8 M^2 n]$ from 
%$[-8M n, 8 M n]^{d-1} \times \{8 Mn\}$ to `the blue part of the boundary of that box' (as indicated in the top half of
%the picture on 18 Aug. (11)), and $F_3$ is the same event `reflected' in the $z_d = 0$ hyperplane.

\noindent
It is easy to see that the event $A_2(8Mn, 8M^2n)$ contains the event $E_1 \cap F_2 \cap F_3$ (see again Figure 1, where the red structure
indicates the above mentioned
`barrier', and the two green paths represent the events $F_2$ and $F_3$), and hence, by independence,

\begin{equation} \label{eq-near-conclusion}
\P\left(A_2(8Mn, 8M^2n)\right) \geq \P(E_1 \cap F_2 \cap F_3) \geq \ep^* \, (\ep')^2.
\end{equation}

Concluding, we have that for each $n \geq 2$, \eqref{l-option1} or \eqref{eq-near-conclusion} holds. From this, 
Proposition \ref{MainProp} follows immediately. \qed

\begin{figure}
	\begin{center}
		\includegraphics[width=0.8\linewidth]{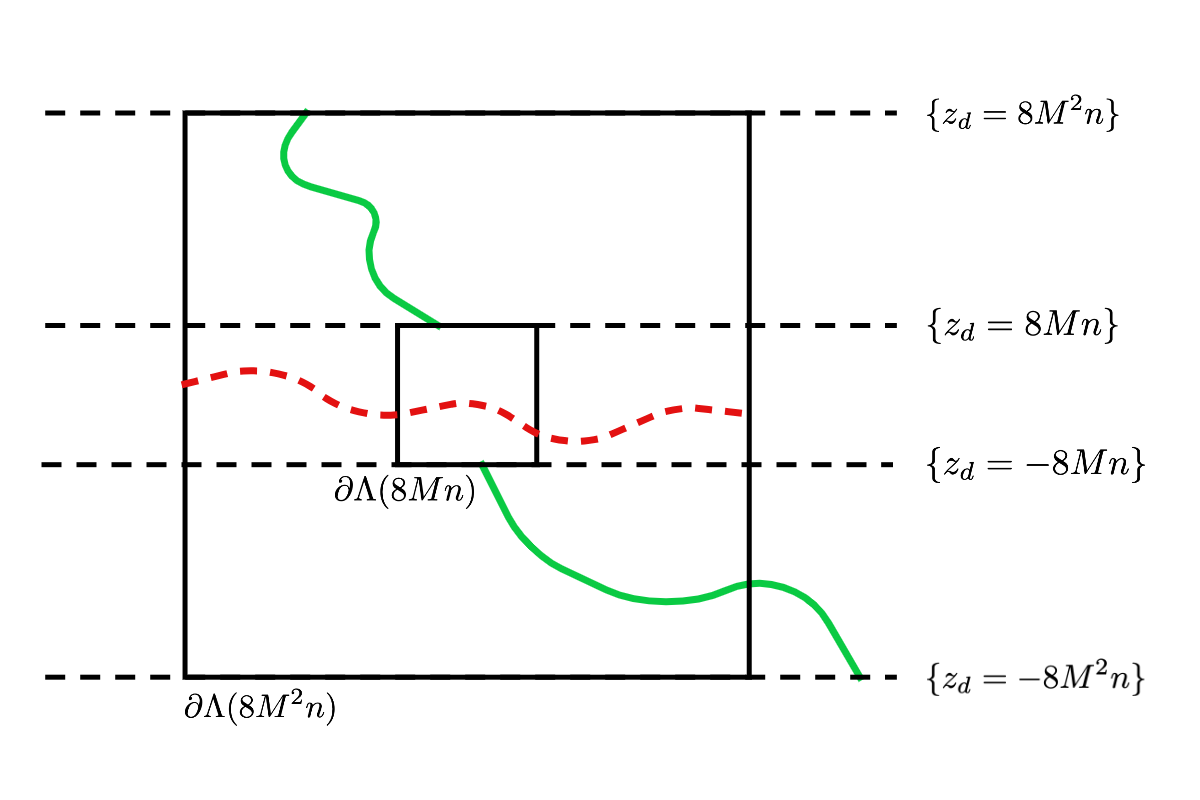}
		\caption{The events $E_1, F_2$ and $F_3$ imply the event $A_2(8Mn, 8M^2n)$.}
		\label{fig: gluing 1}
	\end{center}
\end{figure}

\end{subsection}

\end{section}

\end{document}